\documentclass{amsart}
\usepackage{amssymb,amscd}
\usepackage[all]{xy}

\newtheorem{theorem}{Theorem}[section]

\newtheorem{proposition}[theorem]{Proposition}
\newtheorem{corollary}[theorem]{Corollary}

\theoremstyle{definition}
\newtheorem{definition}[theorem]{Definition}
\newtheorem{remark}[theorem]{Remark}
\newtheorem{example}[theorem]{Example}
\newtheorem{question}[theorem]{Question}

\numberwithin{equation}{theorem}

\def\ge{\geqslant}
\def\le{\leqslant}
\def\dlim{\varinjlim}
\def\phi{\varphi}
\def\epsilon{\varepsilon}
\def\tilde{\widetilde}
\def\to{\longrightarrow}
\def\mapsto{\longmapsto}
\def\GR{\operatorname{GR}}
\def\t{\theta}

\newcommand{\Proj}{\operatorname{Proj}}

\newcommand{\PP}{{\mathbb P}}
\newcommand{\QQ}{{\mathbb Q}}
\newcommand{\NN}{{\mathbb N}}
\newcommand{\ZZ}{{\mathbb Z}}

\newcommand{\calL}{{\mathcal L}}
\newcommand{\calM}{{\mathcal M}}
\newcommand{\calO}{{\mathcal O}}

\newcommand{\fraka}{{\mathfrak a}}
\newcommand{\frakp}{{\mathfrak p}}
\newcommand{\frakm}{{\mathfrak m}}

\begin{document}


\title{Annihilators of Local Cohomology in Characteristic Zero}

\author{Paul Roberts}
\address{Department of Mathematics, University of Utah, 155~South 1400~East, Salt Lake City, UT~84112, USA} \email{roberts@math.utah.edu}
\author{Anurag K. Singh}
\address{Department of Mathematics, University of Utah, 155~South 1400~East, Salt Lake City, UT~84112, USA} \email{singh@math.utah.edu}
\author{V. Srinivas}
\address{School of Mathematics, Tata Institute of Fundamental Research, Homi Bhabha Road, Mumbai~400005, India} \email{srinivas@math.tifr.res.in}

\thanks{P.R.~and A.K.S.~were supported in part by grants from the National Science Foundation. V.S.~was supported by a Swarnajayanthi Fellowship of the DST}

\dedicatory{To Phil Griffith}

\begin{abstract}
This paper discusses the problem of whether it is possible to annihilate elements of local cohomology modules by elements of arbitrarily small order under a fixed valuation. We first discuss the general problem and its relationship to the Direct Summand Conjecture, and next present two concrete examples where annihilators with small order are shown to exist. We then prove a more general theorem, where the existence of such annihilators is established in some cases using results on abelian varieties and the Albanese map.
\end{abstract}

\maketitle

\section{Almost vanishing of local cohomology}

The concept of almost vanishing that we use here comes out of recent work on \emph{Almost Ring Theory} by Gabber and Ramero \cite{GR}. This theory was developed to give a firm foundation to the results of Faltings on \emph{Almost \'etale extensions} \cite{Faltings}, and these ideas have their origins in a classic work of Tate on \emph{$p$-divisible groups} \cite{Tate}. The use of the general theory, for our purposes, is comparatively straightforward, but it illustrates the main questions in looking at certain homological conjectures, as discussed later in the section. The approach is heavily influenced by Heitmann's proof of the Direct Summand Conjecture for rings of dimension three \cite{Heitmann-dim3}.

Let $A$ be an integral domain, and let $v$ be a valuation on $A$ with values in the abelian group of rational numbers; more precisely, $v$ is a function from $A$ to $\QQ\cup\{\infty\}$ such that
\begin{enumerate}
\item $v(a)=\infty$ if and only if $a=0$,
\item $v(ab)=v(a)+v(b)$ for all $a,b\in A$, and
\item $v(a+b)\ge\min\{v(a),v(b)\}$ for all $a,b\in A$.
\end{enumerate}
We will also assume that $v(a)\ge 0$ for all elements $a\in A$.

\begin{definition}\label{def:almostzero}
An $A$-module $M$ is \emph{almost zero} if for every $m\in M$ and every real number $\epsilon>0$, there exists an element $a$ in $A$ with $v(a)<\epsilon$ and $am=0$. When it is necessary to specify the valuation, we say that \emph{$M$ is almost zero with respect to the valuation $v$}.
\end{definition}

We note some properties of almost zero modules:
\begin{enumerate}
\item For an exact sequence
\[
0\to M'\to M\to M''\to 0\,,
\]
the module $M$ is almost zero if and only if each of $M'$ and $M''$ is almost zero.

\item If $\{M_i\}$ is a directed system consisting of almost zero modules, then its direct limit $\dlim_iM_i$ is almost zero.
\end{enumerate}

In \cite{GR} Gabber and Ramero define a module to be almost zero if it is annihilated by a fixed ideal $\frakm$ of $A$ with $\frakm=\frakm^2$. This set of modules also satisfies conditions (1) and (2), though in many cases their condition is stronger than the one in Definition~\ref{def:almostzero}.

The \emph{absolute integral closure} $R^+$ of a domain $R$ is the integral closure of $R$ in an algebraic closure of its fraction field. An important situation for us will be where $(R,\frakm)$ is a complete local ring. In this case, fix a valuation $v\colon R\to\ZZ\cup\{\infty\}$ which is positive on $\frakm$. By Izumi's Theorem \cite{Izumi}, two such valuations are bounded by constant multiples of each other. Since $R^+$ is an integral extension, $v$ extends to a valuation $v\colon R^+\to\QQ\cup\{\infty\}$. Let $A$ be a subring of $R^+$ containing $R$; we often take $A$ to be $R^+$. Note that $v$ is positive on the maximal ideal of $A$. The ring $A$ need not be Noetherian, and by a \emph{system of parameters} for $A$, we shall mean a system of parameters for some Noetherian subring of $A$ that contains $R$.

The main question we consider is whether the local cohomology modules $H^i_\frakm(A)$ are almost zero for $i<\dim A$. Let $x_1\dots,x_d$ be a system of parameters for $R$. Then the local cohomology module $H^i_\frakm(A)$ is the $i$-th cohomology modules of the \v Cech complex
\[
0\to A\to\oplus A_{x_i}\to\oplus A_{x_ix_j}\to\cdots\to A_{x_1\cdots x_d}\to 0\,.
\]
The question whether $H^i_\frakm(A)$ is almost zero for $i=0,\dots,d-1$ is closely related to the question whether the $x_i$ come close to forming a regular sequence in the following sense.

\begin{definition}
A sequence of elements $x_1,\dots,x_d\in A$ is an \emph{almost regular sequence} if for each $i=1,\dots,d$, the module
\[
((x_1,\dots,x_{i-1}):_Ax_i)/(x_1,\dots,x_{i-1})
\]
is almost zero. If every system of parameters for $A$ is an almost regular sequence, we say that $A$ is \emph{almost Cohen-Macaulay}.
\end{definition}

The usual inductive argument as in \cite[Theorem~IV.2.3]{Serre} shows that if $A$ is almost Cohen-Macaulay, then the modules $H^i_\frakm(A)$ are almost zero for $i<\dim A$. However, we do not know whether the converse holds in general.

As motivation for the definitions introduced above, we discuss how these are related to the homological conjectures. Let $x_1,\dots,x_d$ be a system of parameters for a local ring $R$. Hochster's \emph{Monomial Conjecture} states that
\[
x_1^t\cdots x_d^t\notin\big(x_1^{t+1},\dots,x_d^{t+1}\big)R\qquad\text{ for all }
t\ge 0.
\]
This is known to be true for local rings containing a field, and Heitmann \cite{Heitmann-dim3} proved it for local rings of mixed characteristic of dimension up to three. It remains open for mixed characteristic rings of higher dimension, where it is equivalent to several other conjectures such as the Direct Summand Conjecture (which states that regular local rings are direct summands of their module-finite extension rings), the Canonical Element Conjecture, and the Improved New Intersection Conjecture; for some of the related work, we mention \cite{EG, Hochster-CBMS, Hochster-ds, PS} and \cite{Roberts-intersection}.

The connection between the Monomial Conjecture and the almost Cohen-Macaulay property is evident from the following proposition.

\begin{proposition}
Let $R$ be a local domain with an integral extension which is almost Cohen-Macaulay. Then the Monomial Conjecture holds for $R$, i.e., for each system of parameters $x_1,\dots,x_d$ of $R$, we have
\[
x_1^t\cdots x_d^t\notin\big(x_1^{t+1},\dots,x_d^{t+1}\big)R\qquad\text{ for all }t\ge 0\,.
\]
\end{proposition}

\begin{proof}
Let $A$ be an integral extension of $R$ which is almost Cohen-Macaulay with respect to a valuation $v$ which is positive on the maximal ideal of $R$. Then $v(x_i)>0$ for each $i=1,\dots,d$; let $\epsilon$ be the minimum of these positive rational numbers. If $x_1^t\cdots x_d^t\in(x_1^{t+1},\dots, x_d^{t+1})R$ for some integer $t$, then
\[
x_1^t\cdots x_d^t = a_1x_1^{t+1}+\cdots+a_dx_d^{t+1}
\]
for elements $a_i$ of $A$. (The $a_i$ can be chosen in $R$, though we will only consider them as elements of $A$.) Rearranging terms in the above equation, we have
\[
x_1^t(x_2^t\cdots x_d^t-a_1x_1)\in\big(x_2^{t+1},\dots,x_d^{t+1}\big)A\,.
\]
Since $A$ is almost Cohen-Macaulay, the elements $x_1^t,x_2^{t+1},\dots,x_d^{t+1}$ form an almost regular sequence. Hence there exists $c_1\in A$ with $v(c_1)<\epsilon/d$ and
\[
c_1(x_2^t\cdots x_d^t-a_1x_1)\in\big(x_2^{t+1},\dots,x_d^{t+1}\big)A\,.
\]
This implies that $c_1x_2^t\cdots x_d^t\in(x_1,x_2^{t+1},\dots,x_d^{t+1})A$. We now repeat the process for $x_2$, i.e., we have
\[
c_1x_2^t\cdots x_d^t=b_1x_1+b_2x_2^{t+1}+\cdots+b_dx_d^{t+1}
\]
with $b_i\in A$, so
\[
x_2^t(c_1x_3^t\cdots x_d^t-b_2x_2)\in\big(x_1,x_3^{t+1},\dots,x_d^{t+1}\big)A\,.
\]
By an argument similar to the one above, there is an element $c_2\in A$ with $v(c_2)<\epsilon/d$ and
\[
c_1c_2x_3^t\cdots x_d^t\in\big(x_1,x_2,x_3^{t+1},\dots,x_d^{t+1}\big)A\,.
\]
Repeating this procedure $d-2$ more times, we obtain elements $c_1,c_2,\dots,c_d$ in $A$ with $v(c_i)<\epsilon/d$ and
\[
c_1c_2\cdots c_d=u_1x_1+\cdots+u_dx_d
\]
for some $u_i\in A$. But then
\[
v(c_1\cdots c_d)=v(c_1)+\cdots+v(c_d)<d(\epsilon/d)=\epsilon
\]
whereas, since $v(u)\ge 0$ for all $u\in A$, we also have
\[
v(u_1x_1+\cdots+u_dx_d)\ge\min\{v(u_ix_i)\}\ge\min\{v(u_i)+v(x_i)\}
\ge\min\{v(x_i)\}=\epsilon\,,
\]
which is a contradiction.
\end{proof}

To put the results of the remainder of the paper in context, we recall the situation in positive characteristic. Let $R$ be a complete local domain containing a field of characteristic $p>0$. We let $R_\infty$ denote the perfect closure of $R$, that is, $R_\infty$ is the ring obtained by adjoining to $R$ the $p^n$-th roots of all its elements.

\begin{proposition}\label{prop:Rinfinity}
Let $(R,\frakm)$ be a complete local domain containing a field of prime characteristic. Then $R_\infty$ is almost Cohen-Macaulay with respect to any valuation which is positive on $\frakm$.
\end{proposition}

\begin{proof}
Let $v$ be such a valuation, and let $x_1,\dots,x_d$ be a system of parameters for $R_\infty$. Suppose that
\begin{equation}\label{eqn:Rinfinity}
ax_i=b_1x_1+\cdots+b_{i-1}x_{i-1}
\end{equation}
for $a,b_j\in R_\infty$. Let $R'$ be the Noetherian subring of $R_\infty$ generated over $R$ by $a$, $b_1,\dots,b_{i-1}$, and $x_1,\dots,x_d$. By Cohen's structure theorem, $R'$ is a finite extension of a power series ring $S=K[[x_1,\dots,x_d]]$, where $K$ is a coefficient field. Let $m$ be the largest integer such that $R'$ contains a free $S$-module of rank $m$. In this case, the cokernel of
\[
S^m\subseteq R'
\]
is a torsion $S$-module, so there exists a nonzero element $c\in S$ such that $cR'\subseteq S^m$. Taking $p^n$-th powers in equation~\eqref{eqn:Rinfinity} gives us
\[
a^{p^n}x_i^{p^n}\in\big(x_1^{p^n},\dots,x_{i-1}^{p^n}\big)R'\qquad\text{ for all }n\ge0\,.
\]
Multiplying the above by $c$ and using $cR'\subseteq S^m$, we get
\[
ca^{p^n}x_i^{p^n}\in\big(x_1^{p^n},\dots,x_{i-1}^{p^n}\big)S^m\,.
\]
Since $x_1^{p^n},\dots,x_i^{p^n}$ is a regular sequence on the free module $S^m$, it follows that
\[
ca^{p^n}\in\big(x_1^{p^n},\dots,x_{i-1}^{p^n}\big)S^m\subseteq\big(x_1^{p^n},\dots,x_{i-1}^{p^n}\big)R'\,.
\]
Taking $p^n$-th roots in an equation for $ca^{p^n}\in(x_1^{p^n},\dots,x_{i-1}^{p^n})R'$ gives us
\[
c^{1/p^n}a\in\big(x_1,\dots,x_{i-1}\big)R_\infty\qquad\text{ for all }n\ge0\,.
\]
Since the limit of $v(c^{1/p^n})$ is zero as $n\to\infty$, it follows that $R_\infty$ is almost Cohen-Macaulay.
\end{proof}

In \cite{HHbig} Hochster and Huneke proved the much deeper fact that for an excellent local domain $R$ of positive characteristic, the ring $R^+$ is Cohen-Macaulay; see also \cite{HL}. We remark that the subring $R_\infty$ may not be Cohen-Macaulay in general: if $R$ is an $F$-pure ring which is not Cohen-Macaulay, then, since $R\hookrightarrow R_\infty$ is pure, $R_\infty$ is not Cohen-Macaulay as well.

If $R$ is a local domain containing a field of characteristic zero, then $R^+$ is typically not a big Cohen-Macaulay algebra. For example, let $R$ be a normal ring of characteristic zero which is not Cohen-Macaulay. Then the field trace map shows that $R$ splits from finite integral extensions. Consequently a nontrivial relation on a system of parameters for $R$ remains nontrivial in finite extensions, and hence in $R^+$. Specifically, for a ring $(R,\frakm)$ of characteristic zero, the map
\[
H^i_\frakm(R)\to H^i_\frakm(R^+)
\]
is injective for all $i$. This leads to the following question.

\begin{question}\label{q1}
Let $(R,\frakm)$ be a complete local domain. For $i<\dim R$, is the image of natural map
\[
H^i_\frakm(R)\to H^i_\frakm(R^+)
\]
almost zero?
\end{question}

The answer is affirmative if the ring $R$ contains a field of positive characteristic: this follows from Proposition~\ref{prop:Rinfinity}, or from either of the stronger statements \cite[Theorem~1.1]{HHbig} or \cite[Theorem~2.1]{HL}. If $R$ is a three-dimensional ring of mixed characteristic $p$, Heitmann \cite{Heitmann-dim3} proved that the image of $H^2_\frakm(R)$ in $H^2_\frakm(R^+)$ is killed by $p^{1/n}$ for all integers $n\ge0$; more recently, he proved the stronger statement \cite[Theorem~2.9]{Heitmann2005} that $H^2_\frakm(R^+)$ is annihilated by $c^{1/n}$ for all $c\in\frakm$ and $n\ge0$.
Hence the answer to Question~\ref{q1} is also affirmative for mixed characteristic rings of dimension less than or equal to three. 

\section{Examples}

In this section, we present two nontrivial examples where local cohomology modules of characteristic zero rings are annihilated by elements of arbitrarily small positive order. The examples are $\NN$-graded, and in such cases it is natural to use the valuation arising from the grading: $v(r)$ is the least integer $n$ such that the $n$-th degree component of $r$ is nonzero.

\begin{proposition}\label{prop:normalization}
Let $K$ be a field of characteristic zero, and consider the hypersurface $S=K[x,y,z,w]/(xy-zw)$. For distinct elements $\alpha_i$ of $K$, let $\eta$ be a square root of
\[
\prod_{i=1}^4(x-\alpha_iz)\,.
\]
Then the integral closure of $S[\eta]$ in its field of fractions is the ring
\[
R=S\left[\eta,\frac{w}{x}\eta,\frac{w^2}{x^2}\eta\right]\,.
\]
\end{proposition}

\begin{proof}
The element $(w^2/x^2)\eta$ is integral over $S[\eta]$ since
\[
\left(\frac{w^2}{x^2}\eta\right)^2=\prod_{i=1}^4(w-\alpha_iy)\,.
\]
A similar computation shows that $(w/x)\eta$ is integral over $S[\eta]$, and it remains to prove that the integral closure of $S[\eta]$ is generated by these elements. An element of the fraction field of $S[\eta]$ can be written as $a+b\eta$, with $a$ and $b$ from the fraction field of $S$. Now $a+b\eta$ is integral over $S$ if and only if its trace and norm of are elements of $S$. Since $2$ is a unit in $S$ this is equivalent to $a\in S$ and $b^2\eta^2\in S$. Thus the integral closure of $S[\eta]$ is $S\oplus I\eta$, where $I$ is the fractional ideal consisting of elements $b$ with $b^2\eta^2\in S$.

Since $S$ is a normal domain, $b^2\eta^2$ belongs to $S$ if and only if $v_\frakp(b^2\eta^2)\ge0$ for all valuations $v_\frakp$ corresponding to height one prime ideals $\frakp$ of $S$. Note that $v_\frakp(\eta^2)>0$ precisely for the primes $\frakp_0=(x,z)$ and $\frakp_i=(x-\alpha_iz,w-\alpha_iy)$ for $1\le i\le4$. Since $v_{\frakp_0}(\eta^2)=4$ and $v_{\frakp_i}(\eta^2)=1$ for $1\le i\le4$,
the condition for $b$ to be an element of $I$ is that
\[
v_{\frakp_0}(b)\ge-2\qquad\text{ and }\qquad v_\frakp(b)\ge0\text{ for all }\frakp\neq\frakp_0\,.
\]
This implies that $v_\frakp(bx^2)\ge0$ for all height one primes $\frakp$, i.e., that $bx^2\in S$. Let $b=s/x^2$. Then $v_{(x,w)}(b)\ge0$ implies that $s$ must be in the ideal $(x,w)^2$. Hence $I$ is generated over $S$ by $1$, $w/x$, and $w^2/x^2$.
\end{proof}

\begin{example}\label{ex1}
We continue in the notation of Proposition~\ref{prop:normalization}, i.e., $R$ is the normalization of $S[\eta]$. The ring $R$ is normal by construction, and has dimension three. It follows that $H^0_\frakm(R)=0=H^1_\frakm(R)$, where $\frakm$ is the homogeneous maximal ideal of $R$. We show that there are elements of $R^+$ of arbitrarily small positive order annihilating the image of $H^2_\frakm(R)$ in $H^2_\frakm(R^+)$.

Note that $x$, $y$, $z+w$ form a homogeneous system of parameters for the hypersurface $S$, and hence also for $R$. In the ring $R$, we have a relation on these elements given by the equation
\[
\frac{w}{x}\eta\cdot(z+w)=\eta\cdot y+\frac{w^2}{x^2}\eta\cdot x\,.
\]
This is a nontrivial relation in the sense that $(w/x)\eta$ does not belong to the ideal generated by $x$ and $y$, so the ring $R$ is not Cohen-Macaulay. Consider the element of $H^2_\frakm(R)$ given by this relation; it turns out that $H^2_\frakm(R)$ is a one-dimensional $K$-vector space generated by this element, see Remark~\ref{rem:Segre}.

Let $v$ be the valuation defined by the grading on $R$, i.e., $v$ takes value $1$ on $x$, $y$, $z$, and $w$, and $v(\eta)=2$. We construct elements $x_n$ in finite extensions $R_n$ of $R$ with $v(x_n)=1/2^n$ and $x_n(w/x)\eta\in(x,y)R_n$; it then follows that each $x_n$ annihilates the image of the map $H^2_\frakm(R)\to H^2_\frakm(R^+)$.

Let $R_1$ be the extension ring of $R$ obtained by adjoining $\sqrt{x-\alpha_iz}$ for $1\le i\le 4$ and normalizing. We claim that the element $x_1=\sqrt{x-\alpha_1z}$ multiplies $(w/x)\eta$ into the ideal $(x,y)R_1$. To see this, note that
\begin{multline*}
x_1\frac{w}{x}\eta=x_1\frac{w}{x}\prod_{i=1}^4\sqrt{x-\alpha_iz}
=(x-\alpha_1z)\frac{w}{x}\prod_{i=2}^4\sqrt{x-\alpha_iz}\\
=x\left(\frac{w}{x}\prod_{i=2}^4\sqrt{x-\alpha_iz}\right)
-y\left(\alpha_1\prod_{i=2}^4\sqrt{x-\alpha_iz}\right)\,.
\end{multline*}
The element $x_1$ has $v(x_1)=1/2$. To find an annihilator $x_2$ with $v(x_2)=1/4$, we first write
\[
x-\alpha_3z=\beta(x-\alpha_1z)-\gamma(x-\alpha_2z)
\]
for suitable $\beta,\gamma\in K$, and then factor as a difference of squares to obtain
\begin{multline*}
x-\alpha_3z\\
=\left(\sqrt{\beta(x-\alpha_1z)}+\sqrt{\gamma(x-\alpha_2z)}\right)
\left(\sqrt{\beta(x-\alpha_1z)}-\sqrt{\gamma(x-\alpha_2z)}\right)\,.
\end{multline*}
We let
\[
x_2=\sqrt{\sqrt{\beta(x-\alpha_1z)}+\sqrt{\gamma(x-\alpha_2z)}}\,,
\]
which is an element with $v(x_2)=1/4$. Now
\[
x_2\sqrt{x-\alpha_3z}=\lambda\left(\sqrt{\beta(x-\alpha_1z)}+\sqrt{\gamma(x-\alpha_2z)}\right)
\]
where
\[
\lambda=\sqrt{\sqrt{\beta(x-\alpha_1z)}-\sqrt{\gamma(x-\alpha_2z)}}\,,
\]
and so
\begin{multline*}
x_2\eta=\lambda(x-\alpha_1z)\sqrt{\beta(x-\alpha_2z)(x-\alpha_4z)}\\
+\lambda(x-\alpha_2z)\sqrt{\gamma(x-\alpha_1z)(x-\alpha_4z)}\,.
\end{multline*}
Using this, we get
\begin{align*}
x_2\frac{w}{x}\eta
&=x\left(\lambda\frac{w}{x}\sqrt{\beta(x-\alpha_2z)(x-\alpha_4z)}\right)
-y\left(\lambda\alpha_1\sqrt{\beta(x-\alpha_2z)(x-\alpha_4z)}\right)\\
&\quad+x\left(\lambda\frac{w}{x}\sqrt{\gamma(x-\alpha_1z)(x-\alpha_4z)}\right)
-y\left(\lambda\alpha_2\sqrt{\gamma(x-\alpha_1z)(x-\alpha_4z)}\right)
\end{align*}
and consequently $x_2(w/x)\eta\in(x,y)R_2$, where $R_2$ is the finite extension of $R$ obtained by adjoining the various roots occurring in the previous equation and normalizing.

We describe briefly the process of constructing $x_n$ for $n\ge 3$. The first step is to write $\sqrt{x-\alpha_4z}$ in terms of $\sqrt{x-\alpha_1z}$ and $\sqrt{x-\alpha_2z}$ as we did for $\sqrt{x-\alpha_3z}$ above. This enables us to write $\sqrt{x-\alpha_3z}\sqrt{x-\alpha_4z}$ as a product of four square roots, each of which is a linear combination of $\sqrt{x-\alpha_1z}$ and $\sqrt{x-\alpha_2z}$. We can now repeat the process used to construct $x_2$, essentially replacing $x$ by $\sqrt{x-\alpha_3z}$ and $z$ by $\sqrt{x-\alpha_4z}$. Finally, we can repeat this process indefinitely, obtaining elements $x_n$ with $v(x_n)=1/2^n$ which annihilate the given element of local cohomology.
\end{example}

\begin{remark}\label{rem:Segre}
The ring $R$ in the previous example can be obtained as a Segre products of rings of lower dimension, and we briefly discuss this point of view.

Let $A$ and $B$ be $\NN$-graded normal rings which are finitely generated over a field $A_0=B_0=K$. Their \emph{Segre product} is the ring
\[
R=A\#B=\bigoplus_{n\ge0}A_n\otimes_KB_n\,,
\]
which inherits a natural grading where $R_n=A_n\otimes_KB_n$. If $K$ is algebraically closed then the tensor product $A\otimes_KB$ is a normal ring, and hence so is its direct summand $A\#B$. If $M$ and $N$ are $\ZZ$-graded modules over $A$ and $B$ respectively, then their Segre product is the $R$-module
\[
M\#N=\bigoplus_{n\in\ZZ}M_n\otimes_KN_n\,.
\]
Using $\frakm$ to denote the homogeneous maximal ideal of $R$, the local cohomology modules $H_\frakm^k(R)$ can be computed using the K\"unneth formula due to Goto and Watanabe, \cite[Theorem 4.1.5]{GW}:
\begin{multline*}
H^k_\frakm(R)=\left(A\#H^k_{\frakm_B}(B)\right)\ \oplus\
\left(H^k_{\frakm_A}(A)\#B\right)\\
\oplus\bigoplus_{i+j=k+1}\left(H^i_{\frakm_A}(A)\#H^j_{\frakm_B}(B)\right)\,.
\end{multline*}
It follows that if $A$ and $B$ have positive dimension, then
\[
\dim(A\#B)=\dim A+\dim B-1\,.
\]

We claim that the ring $R$ in Example~\ref{ex1} is isomorphic to the Segre product $A\#B$, where
\[
A=K[a,b,c]/\big(c^2-\prod_{i=1}^4(a-\alpha_ib)\big)
\]
is a hypersurface with $\deg a=\deg b=1$ and $\deg c=2$, and $B=K[s,t]$ is a standard graded polynomial ring. The map
\begin{align*}
x&\mapsto as\,,&y&\mapsto bt\,,&z&\mapsto bs\,,&w&\mapsto at\,,\\
\eta&\mapsto cs^2\,,&(w/x)\eta&\mapsto cst\,,&(w/x)^2\eta&\mapsto ct^2
\end{align*}
extends to a $K$-algebra homomorphism $\phi\colon R\to A\#B$. This is a surjective homomorphism of integral domains of equal dimension, so it must be an isomorphism. Since $A$ and $B$ are Cohen-Macaulay rings of dimension $2$, the K\"unneth formula for $H^2_\frakm(R)$ reduces to
\[
H^2_\frakm(R)=\left(A\#H^2_{\frakm_B}(B)\right)\ \oplus\
\left(H^2_{\frakm_A}(A)\#B\right)\,.
\]
The module $H^2_{\frakm_B}(B)$ vanishes in nonnegative degrees, which implies that $A\#H^2_{\frakm_B}(B)=0$. The component of $H^2_{\frakm_A}(A)$ in nonnegative degree is the one-dimensional vector space spanned by the degree $0$ element
\[
\left[\frac{c}{ab}\right]\in H^2_{\frakm_A}(A)\,.
\]
Hence $H^2_\frakm(R)$ is the one-dimensional vector space spanned by $[c/ab]\otimes 1$. The search for elements $x_n\in R^+$ of small degree annihilating the image of $H^2_\frakm(R)$ in $H^2_\frakm(R^+)$ is essentially the search for homogeneous elements of $A^+$, of small degree, multiplying $c$ into the ideal $(a,b)A^+$.
\end{remark}

\begin{example}\label{ex2}
Let $K$ be an algebraically closed field of characteristic zero, $\theta\in K$ a primitive cube root of unity, and set
\[
A=K[x,y,z]/\big(\t x^3+\t^2 y^3+z^3\big)\,.
\]
Let $R$ be the Segre product of $A$ and the polynomial ring $K[s,t]$. Then $R$ is a normal ring of dimension $3$, and the elements $sx$, $ty$, $sy+tx$ form a homogeneous system of parameters for $R$. Using the K\"unneth formula as in Remark~\ref{rem:Segre}, the local cohomology module $H^2_\frakm(R)$ is a one-dimensional vector space spanned by an element corresponding to the relation
\[
sztz(sy+tx)=(sz)^2ty+(tz)^2sx\,.
\]
To annihilate this relation by an element of $R^+$ of positive degree $\epsilon\in\QQ$, it suffices to find an element $u\in A^+$ of degree $\epsilon$ such that
\[
uz^2\in(x,y)A^+\,;
\]
indeed if $uz^2=vx+wy$ for homogeneous $v,w\in A^+$ of degree $1+\epsilon$, then
\[
(s^{\epsilon}u)(sztz)=(s^{\epsilon}tv)(sx)+(s^{1+\epsilon}w)(ty)\,,
\]
and $s^{\epsilon}tv$ and $s^{1+\epsilon}w$ are easily seen to be integral over $S$.

We have now reduced our problem to working over the hypersurface $A$, where we are looking for elements $u\in A^+$ of small degree which annihilate
\[
\left[\frac{z^2}{xy}\right]\in H^2_{\frakm_A}(A^+)\,.
\]
Let $A_1$ be the extension of $A$ obtained by adjoining $x_1$, $y_1$, $z_1$, where
\[
x_1^3=\t^{1/3}x+\t^{2/3}y\,,\qquad y_1^3=\t^{1/3}x+\t^{5/3}y\,,\qquad z_1^3
=\t^{1/3}x+\t^{8/3}y\,.
\]
Note that $x$ and $y$ can be written as $K$-linear combinations of $x_1^3$ and $y_1^3$. Moreover,
\begin{multline*}
(x_1y_1z_1)^3=\big(\t^{1/3}x+\t^{2/3}y\big)\big(\t^{1/3}x+\t^{5/3}y\big)\big(\t^{1/3}x+\t^{8/3}y\big)\\
=\t x^3+\t^2y^3=-z^3\,,
\end{multline*}
so $z$ belongs to the $K$-algebra generated by $x_1$, $y_1$, and $z_1$. Now
\begin{multline*}
\t x_1^3+\t^2 y_1^3+z_1^3=\t\big(\t^{1/3}x+\t^{2/3}y\big)+\t^2\big(\t^{1/3}x+\t^{5/3}y\big)+\big(\t^{1/3}x+\t^{8/3}y\big)\\
=\big(\t^{4/3}+\t^{7/3}+\t^{1/3}\big)x+\big(\t^{5/3}+\t^{11/3}+\t^{8/3}\big)y=0\,,
\end{multline*}
which implies that
\[
A_1=K[x_1,y_1,z_1]/\big(\t x_1^3+\t^2 y_1^3+z_1^3\big)
\]
is a ring isomorphic to $A$. Thus $A\subset A_1$ gives a finite embedding of $A$ into itself under which the generators of degree $1$ go to elements of degree $3$; or, in terms of the original degree, the new generators of the homogeneous maximal ideal have degree $1/3$. Since $[H^2_{\frakm_A}(A)]_0$ is annihilated by all elements of positive degree, the image of $[H^2_{\frakm_A}(A)]_0$ in $H^2_{\frakm_A}(A_1)$ is annihilated by elements of degree $1/3$. Iterating this construction, we conclude that there are elements of arbitrarily small positive degree annihilating the image of $[H^2_{\frakm_A}(A)]_0$ in $H^2_{\frakm_A}(A^+)$. Quite explicitly, we have a tower of extensions
\[
A=A_0\subset A_1\subset A_2\subset\dots\quad\text{ where }\quad
A_n=K[x_n,y_n,z_n]/\big(\t x_n^3+\t^2 y_n^3+z_n^3\big)\,.
\]
The maps $H^2_\frakm(A_n)\to H^2_\frakm(A_{n+1})$ preserve degrees, so $[H^2_{\frakm_A}(A)]_0$ maps to the socle of $H^2_\frakm(A_n)$ which is killed by all elements of $A_n$ of positive degree, e.g., by the elements $x_n,y_n,z_n$ which have degree $1/3^n$.
\end{example}

\begin{remark} 
In \cite[Theorem~2.9]{Heitmann2005} Heitmann proves that if $(R,\frakm)$ is a mixed characteristic excellent local domain of dimension three, then the image of $H^2_\frakm(R)$ in $H^2_\frakm(R^+)$ is annihilated by arbitrarily small powers of every non-unit. The corresponding statement is false for three-dimensional domains of characteristic zero: for the ring $R$ of Example~\ref{ex2}, we claim that $\sqrt{sz}$ does not annihilate the image of $H^2_\frakm(R)\to H^2_\frakm(R^+)$. Because of the splitting provided by field trace, it suffices to verify that
\[
\sqrt{sz}\left(sztz\right)\notin(sx,ty)T\,,
\]
where $T$ is any normal subring of $R^+$ containing $R[\sqrt{sz}]$. Take $T$ to be the Segre product of $\tilde{A}=A[\sqrt{x},\sqrt{y},\sqrt{z}]$ and $\tilde{B}=B[\sqrt{s},\sqrt{t}]$. Note that $\tilde{B}$ is a polynomial ring in $\sqrt{s}$ and $\sqrt{t}$, and that $\tilde{A}$ is the hypersurface
\[
K[\sqrt{x},\sqrt{y},\sqrt{z}]/\left(\t(\sqrt{x})^6+\t^2(\sqrt{y})^6+(\sqrt{z})^6\right)\,.
\]
It is enough to check that $\sqrt{sz}\left(sztz\right)\notin(sx,ty)(\tilde{A}\otimes_K\tilde{B})$, and after specializing $\sqrt{s}\mapsto 1$ and $\sqrt{t}\mapsto 1$ to check that
\[
(\sqrt{z})^5\notin\left((\sqrt{x})^2,(\sqrt{y})^2\right)\tilde{A}\,,
\]
which is immediately seen to be true. The same argument shows that $\sqrt{sx}$, $\sqrt{sy}$, etc. do not annihilate the image of $H^2_\frakm(R)\to H^2_\frakm(R^+)$. The situation is quite similar with Example~\ref{ex1}.
\end{remark}

\section{Annihilators using the Albanese map}

For an $\NN$-graded domain $R$ which is finitely generated over a field $R_0$, let $R^{+\GR}$ be the $\QQ_{\ge0}$-graded ring generated by those elements of $R^+$ which can be assigned a degree such that they satisfy a homogeneous equation of integral dependence over $R$. If $R_0$ is a field of prime characteristic, Hochster and Huneke \cite[Theorem~5.15]{HHbig} proved that the induced map
\[
H^i_\frakm(R)\to H^i_\frakm(R^{+\GR})
\]
is zero for all $i<\dim R$. Translating to projective varieties, one immediately has the vanishing theorem:

\begin{theorem}\cite[Theorem~1.2]{HHbig}
Let $X$ be an irreducible closed subvariety of $\PP^n_K$, where $K$ is a field of positive characteristic. Then for all integers $i$ with $0<i<\dim X$, and all integers $t$, there exists a projective variety $Y$ over a finite extension of $K$ with a finite surjective morphism $f\colon Y\to X$, such that the induced map
\[
H^i(X,\calO_X(t))\to H^i(Y,f^*\calO_X(t))
\]
is zero.
\end{theorem}

Over fields of characteristic zero, the corresponding statements are false because of the splitting provided by field trace. However, the following graded analogue of Question~\ref{q1} remains open.

\begin{question}\label{q:gr1}
Let $R$ be an $\NN$-graded domain, finitely generated over a field $R_0$ of characteristic zero. For $i<\dim R$, is every element of the image of
\[
H^i_\frakm(R)\to H^i_\frakm(R^{+\GR})
\]
killed by elements of $R^{+\GR}$ of arbitrarily small positive degree?
\end{question}

This question, when considered for Segre products, leads to the following:

\begin{question}\label{q:gr2}
Let $R$ be an $\NN$-graded domain of dimension $d$, finitely generated over a field $R_0$ of characteristic zero. Is the image of
\[
\left[H^d_\frakm(R)\right]_{\ge0}\to H^d_\frakm(R^{+\GR})
\]
killed by elements of $R^{+\GR}$ of arbitrarily small positive degree?
\end{question}

In Examples~\ref{ex1} and~\ref{ex2}, we obtained affirmative answers to Question~\ref{q:gr1} by explicitly constructing the annihilators. In this section, we obtain an affirmative answer for the image of $[H^2_\frakm(R)]_0$ and also settle Question~\ref{q:gr2} for rings of dimension two. We first recall some basic facts about the relationship between graded rings and very ample divisors.

If $R$ is a standard graded ring, the associated projective scheme $X=\Proj R$ has a very ample line bundle $\calO_X(1)$ with sections defined by elements of degree one, which generate the line bundle. Conversely, a very ample line bundle defines a standard graded ring and an embedding of $X$ into projective space. 

The strategy is to find, for an arbitrarily large positive integer $n$, a finite surjective map from an integral projective scheme $Y$ to $X$, together with an ample line bundle $\calL$ on $Y$, such that $\calL^{\otimes n}$ is the pullback of $\calO_X(1)$ and such that a section of $\calL$ annihilates the pullback of the given element of cohomology. This will essentially be accomplished by mapping $X$ to its Albanese variety, and pulling back by the multiplication by $N$ map for large integers $N$. The precise result that we prove is as follows.

\begin{theorem}\label{thm:main}
Let $R$ be an $\NN$-graded domain which is finitely generated over a field $R_0$ of characteristic $0$. Let $X=\Proj R$ and let $\eta$ be an element of $H^1(X,\calO_X)$. Then, for every $\epsilon>0$, there exists a finite extension $R\subseteq S$ of $\QQ$-graded domains such that the image of $\eta$ under the induced map
\[
H^1(X,\calO_X)\to H^1(Y,\calO_Y)\qquad\text{ where }Y=\Proj S
\]
is annihilated by an element of $S$ of degree less that $\epsilon$.
\end{theorem}

We remark that since $H^1(X,\calO_X)$ corresponds to the component of $H_\frakm^2(R)$ of degree zero, this theorem only implies that we can annihilate elements of $H^2_\frakm(R)$ of degree zero by elements of small degree. If $H^2_\frakm(R)$ is generated by its degree zero elements---and this happens in several interesting examples---we can deduce the result for all elements of $H^2_\frakm(R)$.

\begin{proof} Replacing $R$ by its normalization, it suffices to work throughout with normal rings. We also reduce to the case where $R$ is a standard $\NN$-graded ring as follows. Using \cite[Lemme~2.1.6]{EGA2}, $R$ has a Veronese subring $R^{(t)}$ which is generated by elements of equal degree. The local cohomology of $R^{(t)}$ supported at its homogeneous maximal ideal $\frakm$ can be obtained by \cite[Theorem~3.1.1]{GW} which states that
\[
H^i_\frakm(R^{(t)})=\bigoplus_{n\in\ZZ}[H^i_\frakm(R)]_{nt}\,.
\]
In particular, we have
\[
[H^2_\frakm(R^{(t)})]_0=[H^2_\frakm(R)]_0=H^1(X,\calO_X)\,.
\]

If elements of this cohomology group can be annihilated in graded finite extensions of $R^{(t)}$, then the same can be achieved in extensions of $R$. 

We next treat the special case where $\Proj R$ is itself an abelian variety, which we denote $A$. For each integer $N$, let $[N_A]\colon A\to A$ be the morphism corresponding to multiplication by $N$. Assume further that the very ample sheaf $\calO_A(1)$ defining the graded ring $R$ satisfies the condition that

\begin{equation}\label{eq:neg1}
[(-1)_A]^*(\calO_A(1))=\calO_A(1)\,.
\end{equation}
Note that, if ${\mathcal L}$ is any very ample line bundle on $A$, the new very ample line bundle \[\calO_A(1)={\mathcal L}\otimes [(-1)_A]^*{\mathcal L}\] satisfies this further assumption.

We recall two facts about abelian varieties from Mumford \cite{Mumford}:
\begin{enumerate}

\item $H^1(A,\calO_A(1))=0$. By ``The Vanishing Theorem" \cite[page 150]{Mumford}, given a line bundle $\calL$, there is a unique integer $i$ such that $H^i(A,\calL)$ is nonzero. Since $\calO_A(1)$ is very ample, this integer must be $0$.

\item $[N_A]^*(\calO_A(1))=\calO_A(N^2)$. This follows from \cite[Corollary~II.6.3]{Mumford} since we are assuming \eqref{eq:neg1}. 
\end{enumerate}

The theorem in this case follows from these two properties: the morphism $[N_A]$ induces a map
\[
R=\bigoplus_n\Gamma(A,\calO_A(n))\to\bigoplus_n\Gamma(A,[N_A]^*\calO_A(n))\,,
\]
and, by the second property above,
\[
\Gamma(A,[N_A]^*\calO_A(n))=\Gamma(A,\calO_A(N^2n))\,.
\]

Thus we have a map of graded rings from $R$ to itself, that takes an element of degree $1$ to an element of degree $N^2$. Denote the new copy of $R$ by $S$, and regrade $S$ with a $\QQ$-grading such that the map $R\to S$ preserves degrees. This implies that $S$ has elements $s$ of degree $1/N^2$ under the new grading. Such an element $s$ must annihilate the image of $\eta\in H^1(A,\calO_A)$, since the product $s\cdot\eta$ lies in $H^1(A,\calO_A(1))=0$. Hence for each positive integer $N$, we have found a finite extension of $R$ with an element of degree $1/N^2$ that annihilates the image of $\eta$.

The remainder of the proof is devoted to reducing to the previous case. Let $R$ be a graded domain such that $X=\Proj R$ is normal. Let $A$ be the \emph{strict Albanese variety} of $X$ as defined in \cite{Chevalley}. It is the dual abelian variety to the Picard variety of $X$ (in the sense of Chevalley-Grothendieck) which parametrizes line bundles algebraically equivalent to 0. Let $\phi\colon X\to A$ be the corresponding Albanese morphism. Then (since the ground field has characteristic 0) $\phi$ induces an isomorphism
\[
H^1(A,\calO_A)\cong H^1(X,\calO_X)\,,
\]
see Chevalley \cite{Chevalley}. 

Since $A$ is an abelian variety, it has a very ample invertible sheaf $\calO_A(1)$, see for example \cite[pp.~60--62]{Mumford}. After replacing $\calO_A(1)$ by $\calO_A(1)\otimes [(-1)_A]^*\calO_A(1)$ if necessary, we may assume that
\[
[(-1)_A]^*(\calO_A(1))\cong \calO_A(1)\,.
\]
We let $\calO_X(1)$ denote the very ample invertible sheaf defined by the grading on $R$. Let $\pi_1\colon Y_1\to X$ be the pullback of multiplication by $N$ on $A$, and let $\phi_1\colon Y_1\to A$ be the map induced by $\phi$, so that we have the fiber product diagram below.
\[
\xymatrix{Y_1\ar[r]^{\phi_1}\ar[d]_{\pi_1}& A\ar[d]^{[N_A]}\\
X\ar[r]^{\phi} & A}
\]

Let $\calM_1=\phi_1^*(\calO_A(1))$. Then
\[
\pi_1^*(\phi^*(\calO_A(1)))=\phi_1^*([N_A]^*(\calO_A(1)))=\phi_1^*(\calO_A(N^2))=\calM_1^{\otimes N^2}.
\]

Now let $m$ be an integer such that $\phi^*(\calO_A(-1))\otimes \calO_X(m)$ is globally generated; such an $m$ exists since $\calO_X(1)$ is ample, and we fix one such $m$. 

Since the sheaf $\phi^*(\calO_A(-1))\otimes \calO_X(m)$ is globally generated, there exists a map $\psi\colon X\to\PP^n$ such that
\[
\psi^*(\calO_{\PP^n}(1))=\phi^*(\calO_A(-1))\otimes \calO_X(m)\,.
\]

Let $\alpha\colon\PP^n\to\PP^n$ be a finite map such that $\alpha^*(\calO_{\PP^n}(1))=\calO_{\PP^n}(N)$; for example, we can take $\alpha$ to be the map defined by the ring homomorphism on a polynomial 
ring that sends the variables to their $N$-th powers. 

Let $Y_2$ be the fiber product of $\psi$ and $\alpha$, which gives us a diagram
\[
\xymatrix{Y_2 \ar[r]^{\phi_2}\ar[d]_{\pi_2}& \PP^n\ar[d]^{\alpha} \\
X \ar[r]^{\psi} & \PP^n\,.} 
\]
Let $\calM_2=\phi_2^*(\calO_{\PP^n}(1))$. We then have
\begin{multline*}
\pi_2^*(\phi^*(\calO_A(-1))\otimes\calO_X(m))\cong\pi_2^*(\psi^*(\calO_{\PP^n}(1)))=\phi_2^*(\alpha^*(\calO_{\PP^n}(1)))\\
\cong \phi_2^*(\calO_{\PP^n}(N))\cong \calM_2^{\otimes N}\,.
\end{multline*}
Note that the above morphisms $\pi_i\colon Y_i\to X$, $i=1,2$ are \emph{finite and surjective}. 

Let $Y$ be a component of the normalization of the reduced fibre product of $\pi_1\colon Y_1\to X$ and $\pi_2\colon Y_2\to X$, with the property that $Y\to X$ is surjective. Since $Y_1\times_XY_2\to X$ is finite and surjective, some irreducible component of this fiber product maps onto $X$ via a finite morphism, and we may take $Y$ to be the normalization of any such component. 

We then have an induced finite surjective map $\pi\colon Y\to X$, and induced maps $\mu_1\colon Y\to Y_1$ and $\mu_2\colon Y\to Y_2$, giving a commutative diagram
\[
\xymatrix{Y \ar[r]^{\mu_1}\ar[d]_{\mu_2}\ar[dr]^{\pi}& Y_1\ar[d]^{\pi_1} \\
Y_2 \ar[r]^{\pi_2} & X\,.} 
\]

By construction, we have
\begin{multline*}
\pi^*\calO_X(m)=\pi^*(\phi^*\calO_A(1)\otimes\phi^*\calO_A(-1)\otimes\calO_X(m))\\
=\mu_1^*\pi_1^*(\phi^*(\calO_A(1)))\otimes \mu_2^*\pi_2^*(\phi^*(\calO_A(-1))\otimes\calO_X(m))\\
=\mu_1^*(\calM_1^{\otimes N^2})\otimes \mu_2^*(\calM_2^{\otimes N}) = \calM^{\otimes N}\,,
\end{multline*}
where $\calM=\mu_1^*\calM_1^{\otimes N}\otimes \mu_2^*\calM_2.$

Now $\calM_1= \phi_1^*(\calO_A(1))$ is generated by global sections of the form $\phi_1^*(u)$ with $u\in H^0(A,\calO_A(1))$. Choose any such element $u$ such that its image in $H^0(Y,\mu_1^*\calM_1)$ is nonzero. Let $v$ be a nonzero element of $H^0(Y,\mu_2^*\calM_2)$, and let $s$ be the image of $\mu_1^*\phi_1^*(u^N)\otimes v$ in $H^0(Y,\calM)$. Then $s$ is a nonzero section of $\calM$, and since $\pi^*(\calO_X(m))=\calM^{\otimes N}$, the degree of $s$ in the grading induced from that on $R$ is $m/N$. We claim that the composition
\[
\CD
H^1(X,\calO_X)@>\pi^*>>H^1(Y,\calO_Y)@>\cdot s>>H^1(Y,\calM)
\endCD
\]
vanishes. For this, it suffices to show that the composition
\[
\CD
H^1(X,\calO_X)@>\pi^*>>H^1(Y,\calO_Y)@>\mu_1^*\phi_1^*(u)>>H^1(Y,\mu_1^*\calM_1)
\endCD
\]
vanishes, which in turn reduces to showing that
\[
\CD
H^1(X,\calO_X)@>\pi_1^*>>H^1(Y_1,\calO_{Y_1})@>\phi_1^*(u)>>H^1(Y_1,\calM_1)
\endCD
\]
vanishes. Since $\phi^*\colon H^1(A,\calO_A)\to H^1(X,\calO_X)$ is an isomorphism, we further reduce to showing that
\[
\CD
H^1(A,\calO_A)@>[N_A]^*>>H^1(A,\calO_A)@>\cdot u>>H^1(A,\calO_A(1))
\endCD
\]
vanishes. But this is true since $H^1(A,\calO_A(1))=0$.

Since we can make $m/N$ arbitrarily small by choosing $N$ large, this completes the proof.
\end{proof}

As a corollary, we see that the answer to Question~\ref{q:gr2} is affirmative for rings of dimension two:

\begin{corollary}\label{cor:dim2}
Let $R$ be an $\NN$-graded domain of dimension $2$, which is finitely generated over a field $R_0$ of characteristic zero. Then the image of
\[
\left[H^2_\frakm(R)\right]_{\ge0}\to H^2_\frakm(R^{+\GR})
\]
is killed by elements of $R^{+\GR}$ of arbitrarily small positive degree.
\end{corollary}

\begin{proof}
By adjoining roots of elements if necessary, we may assume that $R$ has a system of parameters $x,y$ consisting of linear forms. Theorem~\ref{thm:main} implies that the image of $[H^2_\frakm(R)]_0$ is killed by elements of $R^{+\GR}$ of arbitrarily small positive degree, so it suffices to prove that $[H^2_\frakm(R)]_{\ge0}$ is the $R$-module generated by $[H^2_\frakm(R)]_0$.

Since $x$ and $y$ are linear forms, we have $[R_{xy}]_{n+1}=[R_{xy}]_n\cdot R_1$ for all integers $n$. Computing $H^2_\frakm(R)$ using the \v Cech complex on $x,y$, it follows that
\[
\left[H^2_\frakm(R)\right]_{n+1}=\left[H^2_\frakm(R)\right]_n\cdot R_1
\qquad\text{ for all }n\in\ZZ\,.\qedhere
\]
\end{proof}

\section{Closure operations}

The issues discussed here are closely related to closure operations considered by Hochster and Huneke, and by Heitmann. The \emph{plus closure} of an ideal $\fraka$ of a domain $R$ is defined as $\fraka^+=\fraka R^+\cap R$. It has desirable properties for rings of prime characteristic, e.g., it bounds colon ideals on systems of parameters: if $x_1,\dots,x_d$ is a system of parameters for an excellent local domain $R$ containing a field of prime characteristic, then
\[
(x_1,\dots,x_{i-1}):_Rx_i\subseteq (x_1,\dots,x_{i-1})^+\qquad\text{ for all }i\,.
\]
In general, plus closure does not have this colon-capturing property for rings of characteristic zero or of mixed characteristic. Several alternative closure operations are defined by Heitmann in \cite{Heitmann2001} including the \emph{extended plus closure}. Building on these ideas, he settled the Direct Summand Conjecture for mixed characteristic rings of dimension three \cite{Heitmann-dim3}. In \cite[Theorem~1.3]{Heitmann2005} Heitmann proved that extended plus closure has the colon-capturing property for arbitrary sets of three parameters in excellent domains of mixed characteristic.

Let $(R,\frakm)$ be a complete local domain and fix, as usual, a valuation $v\colon R\to\ZZ\cup\{\infty\}$ which is positive on $\frakm$ and extend to $v\colon R^+\to\QQ\cup\{\infty\}$. In \cite{HHJPAA} Hochster and Huneke define the \emph{dagger closure} $\fraka^\dagger$ of an ideal $\fraka$ as the ideal consisting of all elements $x\in R$ for which there exist elements $u\in R^+$, of arbitrarily small positive order, with $ux\in\fraka R^+$. In \cite[Theorem~3.1]{HHJPAA} it is proved that the dagger closure $\fraka^\dagger$ agrees with the tight closure $\fraka^*$ for ideals of complete local domains of prime characteristic, see also \cite[\S\,6]{HHJAMS}. While tight closure is defined in characteristic zero by reduction to prime characteristic, the definition of dagger closure is characteristic-free.

\begin{question}[Hochster-Huneke]\label{q:dagger}
Does the dagger closure operation have the colon-capturing property, i.e., if $x_1,\dots,x_d$ is a system of parameters for a complete local domain $R$, is it true that
\[
(x_1,\dots,x_{i-1}):_Rx_i\subseteq(x_1,\dots,x_{i-1})^\dagger\,?
\]
\end{question}

According to Hochster and Huneke \cite[page~244]{HHJPAA} ``it is important to raise (and answer) this question.'' If Question~\ref{q:dagger} has an affirmative answer, then so does Question~\ref{q1}.

Consider the hypersurface $K[[x,y,z]]/(x^3+y^3+z^3)$. If $K$ has prime characteristic, a straightforward calculation---performed in many an introductory lecture on tight closure theory---shows that $z^2\in(x,y)^*$. If $K$ has characteristic zero, the ``reduction modulo $p$'' nature of the definition of tight closure \cite[\S\,3]{HHJAMS} immediately yields $z^2\in (x,y)^*$ once again. In contrast, the computation that $z^2\in (x,y)^\dagger$ is quite delicate and, aside from the linear change of variables, is the computation we performed in Example~\ref{ex2}. While concrete descriptions of the multipliers of small order are available in this and some other examples, dagger closure remains quite mysterious even in simple examples such as diagonal hypersurfaces:

\begin{question}
Let $K$ be a field of characteristic zero, and let
\[
R=K[[x_0,\dots,x_d]]/(x_0^n+\cdots+x_d^n), \qquad\text{where $n>d$}.
\]
Does $x_0^d$ belong to the dagger closure of the ideal $(x_1,\dots,x_d)$?
\end{question}

A routine computation of tight closure shows that $x_0^d\in(x_1,\dots,x_d)^*$. In the case $d=2$, we have $x_0^2\in(x_1,x_2)^\dagger$ by Corollary~\ref{cor:dim2}.


\end{document}